\documentclass[12pt]{article}
\usepackage{theos}
\usepackage{nzqetc}
\newcommand{\tens}{{\otimes}}
\newcommand{\p}{{\mathfrak{p}}}
\newcommand{\m}{{\mathfrak{m}}}
\newcommand{\f}{{\phi}}
\newcommand{\mpr}{{{\m}^{'}}}
\newcommand{\apr}{A^{'}}
\newcommand{\rpr}{R^{'}}
\newcommand{\Rpr}{\rpr_{\mpr}}
\newcommand{\ai}{{\mathfrak{a}}_i}

\begin{document}
\title{On ring homomorphisms of Azumaya algebras}
\author{Kossivi Adjamagbo, Jean-Yves Charbonnel and Arno van den Essen}
\maketitle

\vspace*{1cm}

{\sl {\bf Abstract:} The main theorem (Theorem 4.1) of this paper claims that any ring 
morphism from an Azumaya algebra of constant rank over a commutative ring 
to another one of the same constant rank and over a reduced commutative 
ring induces a ring morphism between the centers of these algebras. The 
second main theorem (Theorem 5.3) implies (through its Corollary 5.4) that 
a ring morphism between two Azumaya algebras of the same constant rank is 
an isomorphism if and only if it induces an isomorphism between their 
centers. As preliminary to these theorems we also prove a theorem (Theorem 
2.8) describing the Brauer group of a commutative 
Artin ring and give an explicit proof of the converse of the
  Artin-Procesi theorem on Azumaya algebras (Theorem 2.6).
  }

\vspace*{1cm}

\section{Introduction}
	According to A. Grothendieck's point of view on algebraic geometry, as 
explained by J. Dieudonn\'e in [DIE], the essence of algebraic 
geometry is the study of scheme morphisms and not of schemes themselves. 
 From this point of view, which attaches more importance to the structure's 
morphisms than to the structures themselves, the most distinguished property 
of Azumaya algebras is that any algebra endomorphism of an Azumaya algebra 
is an automorphism. This is a well-known fact since 1960 thanks to the 
historical paper [AG] of M. Auslander and O. Goldman, which has clearly been
 the source of inspiration of the famous Dixmier 
conjecture formulated in 1968 [DIX] for the first Weyl algebra over a field of 
characteristic zero. This conjecture in its present general formulation claims 
that any algebra endomorphism of a Weyl algebra of index $n$ over such a
 field is an  automorphism.

	This historical link between Azumaya algebras and Weyl algebras has been 
confirmed since 1973 by a remarkable theorem of P. Revoy in [R] asserting 
that any Weyl algebra of index $n$ over a field of positive characteristic is 
a free Azumaya algebra over its center and that the latter is isomorphic 
as algebra over this field to the algebra of polynomials in $2n$ 
indeterminates over this field. This link is not formal but substantial and 
fecund, for it suggests naturally a strategy and many tactics to investigate 
the Dixmier Conjecture with the help of tools from the theory of Azumaya 
algebras which have already been proved and which are to be proved.

	Given an endomorphism $f$ of a Weyl algebra $A_n(k)$ over a field $k$ of 
characteristic zero, this strategy consists in reducing the problem first to an 
endomorphism $f_R$ of a Weyl algebra $A_n(R)$ over a suitable finitely 
generated subring $R$ of $k$, then to an endomorphism $f_p$ of a Weyl algebra 
$A_n(R_p)$ over a suitable reduction $R_p$ modulo a suitable prime integer $p$ of 
the ring $R$ in such a way that $R_p$ is a field of characteristic $p$. From 
then, tools of the theory of Azumaya algebras applied to the Azumaya 
algebra $A_n(R_p)$ suggest many tactics to prove that $f_p$ is an automorphism. 
It remains only to find a good tactic to conclude that $f$ is an automorphism.

	This strategy has been partially applied since 2002 by Y. Tsuchimoto in 
[T], but without the clarity and the power of the tools of Azumaya algebra 
theory.
% and with some lack of rigor like in the proof of the lemma 2 of his 
%paper as he himself recognized it in a footnote with the promise to correct it 
%in a forthcoming paper.

The aim of the present paper is to prove some new tools of Azumaya algebra 
theory needed for the investigation of the Dixmier Conjecture by the evocated 
strategy, and more generally for the deepening, the completeness and the 
beauty of this theory which has already benefited from the contributions of 
an impressive army of researchers among which A. Grothendieck, M. Artin and O. 
Gabber, as explicited in [V]. As indicated in the title of the present paper,
these new tools concern ring morphisms of Azumaya algebras, unlike the cited
 fundamental theorem on algebra endomorphisms of these algebras.

	The main theorem (Theorem 4.1) of this paper claims that any ring 
morphism from an Azumaya algebra of constant rank over a commutative ring 
to another one of the same constant rank and over a reduced commutative 
ring induces a ring morphism between the centers of these algebras. The 
second main theorem (Theorem 5.3) implies (through its Corollary 5.4) that 
a ring morphism between two Azumaya algebras of the same constant rank is 
an isomorphism if and only if it induces an isomorphism between their 
centers. As preliminary to these theorems we also prove a theorem (Theorem 
2.8) describing the Brauer group of a commutative 
Artin ring which implies the fact, well-known 
to the specialist, that the Brauer group of a finite ring is trivial, hence 
that an Azumaya algebra over such a ring is isomorphic to an algebra of 
matrices over this ring (Corollary 2.9). As another preliminary to these
 theorems, we also give an explicit proof of the converse of the
  Artin-Procesi theorem on Azumaya algebras (Theorem 2.6).
 Using this result we sketch an alternative proof of the main theorem. 

	Before getting to the heart of the matter, we would like to express
 our deep gratitude to M. Kontsevich for a fruitful and hearty discussion
  which led us to guess and prove the main theorem below, inspired by
   his conviction of the truth of this statement in the case of 
    matrix rings over fields. We also insist to thank our colleague
  Alberto Arabia from the University Paris 7 for 
a crucial idea in the proof of Lemma 3.2  below.

\section{Preliminaries on Azumaya algebras.}

First we recall some well-known facts about Azumaya
 algebras. Then we add some new ones. For details we refer to [AG], [DI] and
 [KO].

\begin{defn}
Let $R$ be a commutative ring and $A$ an $R$-algebra given by
a ring homomorphism $\phi:R\rightarrow A$ ($A$ need not be commutative).
Then $A$ is called a central $R$-algebra if $A$ is a faithful $R$-module i.e.
$Ann_R(A)=0$ and the center of $A$, denoted $Z(A)$, equals $R. 1_A = 
\phi(R)$ , where $1_A$ denotes the unit element of $A$. Finally, $A$ is 
called an Azumaya algebra over $R$ if $A$ is a
finitely generated faithful $R$-module such that for each
$\m\in Max(R)$ $A/\m A$ is a central simple $R/\m$-algebra.
\end{defn}

\begin{exmp}
(1) Since the ring of $n\times n$ matrices over a field is central simple,
it follows easily that all matrix rings $M_n(R)$ are Azumaya
algebras over $R$.

\medskip

\noindent (2) Other examples of Azumaya algebras are given by
the Weyl algebras $A_n(k)$ over a field $k$ of characteristic
$p>0$. Indeed, it is proved in [R] that 
$A_n(k)=k[x_1,\ldots,x_n,y_1,\ldots,y_n]$,
with the relations $[y_i,x_j]=\delta_{ij}$, is an Azumaya algebra and a 
free module over its center  $k[{x_1}^p,\ldots,{x_n}^p,{y_1}^p,\ldots,{y_n}^p]$ which 
is isomorphic to the $k$-algebra of polynomials in $2 n$ indeterminates. 
\end{exmp}

\begin{rem}
The main well-known properties of Azumaya algebras we need are the 
following (notations as above) :
\medskip

\noindent (1) $A$ is a central $R$-algebra and a finitely generated
 projective $R$-module which contains $R$ as a direct summand.

\medskip

\noindent (2) There is a 1-1 correspondence between the ideals $I$ of $R$ and
the two-sided ideals $J$ of $A$ given by $I\rightarrow IA$ and
$J \rightarrow \phi^{-1}(J)$.

\medskip

\noindent (3) Consequently, if $(\ai)_{i\in I}$ is a family of ideals  in $R$,
then $\cap_i (\ai A)=(\cap\ai)A$,
 since $\phi^{-1}(\cap(\ai A))=\cap\phi^{-1}(\ai A)=\cap\ai$.
\medskip

\noindent (4) If $S$ is a commutative $R$-algebra, then $A\tens_R S$ is an 
Azumaya
algebra over $S$.

\medskip

\noindent (5) Hence, for each $\p\in Spec(R)$ , $A\tens_R R_\p$ is a finite 
free $R_\p$-module of rank $r_\p(A)$  depending only on the connected 
component of $\p$ in $Spec(R)$ and called the rank of $A$ at $\p$.

\medskip

\noindent (6) If $A_1$ and $A_2$ are Azumaya algebras over $R$, then
so is $A_1\tens_RA_2$.

\medskip

\noindent (7) If $A$ is an Azumaya algebra over $R$, then so is $A^0$
and $A\tens_RA^0\simeq End_R(A)$ (here $A^0$ denotes the opposite
algebra of $A$ i.e. $A^0:=A$ as a set, it has the same addition as $A$
and its multiplication is given by $a\star b=b.a$).

\medskip

\noindent (8) If $A_2\subset A_1$ are Azumaya algebras over $R$, then
the commutant of $A_2$ in $A_1$ i.e.

$$A_1^{A_2}:=\{a_1\in A_1\, |\, a_1a_2=a_2a_1 \, for\, all\,  a_2\in A_2\}$$

\noindent is an Azumaya algebra over $R$. Furthermore, the canonical map

\medskip

$A_2\otimes_R A_1^{A_2}\rightarrow A_1$, $a_2\otimes a_1\rightarrow a_1a_1$

\medskip

\noindent is an isomorphism.
\end{rem}

\begin{defn}
\noindent (1) With the notations above, $A$ is said to have constant rank 
$r$ if $r = r_\p(A)$ for each $\p\in Spec(R)$. It is for example the case 
if the reduced ring of $R$ is integral, in particular if $R$ is integral.

\medskip

\noindent (2) On the set of Azumaya algebras one defines an equivalence 
relation by
saying that two Azumaya algebras $A_1$ and $A_2$ are equivalent if
$A_1\tens_R A_2\simeq End_R(P)$ for some faithful finitely generated
projective $R$-module $P$. It then follows from the properties above that
the operation $\tens_R$ defines a group structure on the set of
equivalence classes of Azumaya algebras: the class of $R$ is the neutral
element and the class of $A^0$ is the inverse of the class of $A$. This
abelian group is called the Brauer group of $R$ and is denoted by $Br(R)$.
\end{defn}

\begin{rem}
\noindent (1) Since obviously $A$ is equivalent to $R$ iff $A\simeq 
A\tens_R R\simeq End_R(P)$ for some
faithful finitely generated projective $R$-module $P$, it follows that
$Br(R)$ is trivial iff each Azumaya algebra $A$ over $R$ is isomorphic to
$End_R(P)$ for some faithful finitely generated projective $R$-module $P$.

\medskip

\noindent (2) It follows that if $R$ is a local ring, then $Br(R)$ is 
trivial iff for
each Azumaya algebra $A$ over $R$ there exists an $n\in\nx$ such that
$A\simeq M_n(R)$.

\medskip

\noindent (3) It also follows that if $R$ is a semi-local ring with a 
trivial Brauer group, then for
each Azumaya algebra $A$ of constant rank over $R$ there exists $n\in\nx$ 
such that
$A\simeq M_n(R)$, according to the liberty of projective modules of 
constant rank over semi-local rings (see for instance [BO1], Ch.II,\S 5, no. 3,
 Proposition. 5).

\medskip

\noindent (4) It is proved in [AG], th. 6.5 that the Brauer group of any 
complete local commutative ring is isomorphic to the one of its residue field.

\medskip

\noindent (5) It is easy to check that the Brauer group of the direct sum 
of a finite number of commutative rings is the direct sum of their Brauer 
groups (see for instance [O],p. 58, exercise(i))

\medskip

\noindent (6) One of the deepest results on Azumaya algebras, the
 Artin-Procesi theorem, 
in its version explicitely proved by the authors and their followers, 
claims that if a ring satisfies all the $\zx$-identities of the ring of $n$ 
by $n$ matrices over $\zx$ and if  no nonzero homomorphic image of this ring 
satisfies all the $\zx$-identities of the ring of $n-1$ by $n-1$ matrices 
over $\zx$, then this ring is an Azumaya algebra of constant rank $n^2$ over 
its center (see [AR] and [P]). A recent and short proof can be found in 
[DR], Theorem 9.3.

\medskip

\noindent (7) The following theorem is the converse of this formulation of
the Artin-Procesi theorem (see Theorem 9.5 of [DR]). It has been claimed
 without proof in [AR] and [P] and many followers. It seems that no explicit
  proof has been published. Also the proof given in [DR] is based on unproven 
remarks concerning matrix rings over commutative rings. 
\end{rem}

\begin{thm}
Any Azumaya algebra of constant rank $n^2$ satisfies all the $\zx$-identities 
of the ring $M_n(\zx)$ and no nonzero homomorphic 
image of this ring satisfies all the $\zx$-identities of the ring
$M_{n-1}(\zx)$. 
\end{thm}

\noindent\textbf{Proof}. Let $A$ be such an algebra. According to the 
Amitsur-Levitzki 
theorem ([RO], Theorem 1.4.5 and  Corollary 1.8.43), no nonzero homomorphic
 image of $A$ 
satisfies the standard identity of degree $2n - 2$. Hence, according to the 
Amitsur-Levitzki theorem ([RO], Theorem 1.4.1), no nonzero homomorphic
 image of $A$ 
satisfies all the $\zx$-identities of the ring $M_{n-1}(\zx)$. Finally,
 according to Rowen's version of the Artin-Procesi theorem
 ([RO], Theorem 1.4.1) 
and to the MacConnell-Robson version of the Artin-Procesi theorem ([M], 
Theorem 7.14 and Corollary 6.7), $A$ satisfies all the $\zx$-identities of
 the ring $M_n(R)$, where $R$ is the ring of polynomials with coefficients
  in $\zx$ generated by a 
family of indeterminates indexed by $\{1,...,n\} ^2 \times \nx$. As $R$ is a
 domain, it follows from [M], Theorem 7.14, Corollary 6.7 and 2.2(1)  
that $M_n(R)$ satisfies all the $\zx$-identities of the ring 
$M_n(\zx)$. So we can conclude that $A$ satisfies all the $\zx$-identities of
 the ring $M_n(\zx)$, as desired. 

\begin{exmp}
(See for instance [BO2], Ch.VIII,\S 10, n. 4 or [O], Th. 6.36)

\noindent (1) It follows from two theorems of Wedderburn that the Brauer group 
of any finite field is trivial.

\medskip

\noindent (2) It follows from a theorem of Wedderburn that the Brauer group of 
any algebraically closed field is trivial.

\medskip

%\noindent (3) The Brauer group of the field of real numbers is a cyclique 
%group of order two.

%\medskip

%\noindent (4) The Brauer group of the ring of natural integers is trivial.

%\medskip

%\noindent (5) The Brauer group of the ring of integers of a numbers field 
%with a number $r > 1$ of embeddings into the field of real numbers is the 
%direct sum of $r-1$ cyclic group of order two.
\end{exmp}

\begin{thm}
The Brauer group of a commutative Artin ring is isomorphic to the direct 
sum of the residue fields of its maximal ideals.
\end{thm}

\noindent\textbf{Proof}. Let $R$ be such a ring. According to the structure 
theorem for commutative Artin rings (see for instance [AT], Ch.8, th. 8.7),  
$R$ is the direct sum of complete local rings $R_i$ with the same 
residue  fields $k_i$ as of the maximal ideals of $R$. So according to 2.5(5) 
 the Brauer group of $R$ is the direct sum of the ones of the 
$R_i$'s. But according to 2.5(4) the Brauer group of $R_i$ is 
isomorphic to the one of $k_i$. The conclusion follows.

\begin{cor}
The Brauer group of a finite commutative ring or more generally of a
 commutative Artin ring which is such that each of its residue fields is either
finite or algebraically closed,is trivial. Hence, any Azumaya algebra of
 constant rank over such a ring  is isomorphic to an algebra of matrices
  over this ring.
\end{cor}

\noindent\textbf{Proof}. Follows from 2.8, 2.5(3) and the
examples 2.7(1),(2).

\begin{rem}
\noindent (1) Considering an algebraic closure of the residue field of 
a prime ideal of a commutative ring, it follows from 2.7(2)  
that the rank of any Azumaya algebra at a prime ideal of its center is the 
square of an integer.

\medskip

\noindent (2) To conclude this section let us recall for the purpose of the 
proof of the main theorem below, that a commutative ring is called a 
Jacobson ring if each prime ideal is the
intersection of a family of maximal ideals.

\medskip

\noindent (3) A typical example of such a ring is any finitely generated 
$\zx$-algebra. In this case the residue field of any of its maximal ideals  
is finite (see for instance [BO1], V,\S 3, no. 4, Th.3).
\end{rem}

\section{Ring homomorphisms between matrix rings}

In order to prepare for the proof of the main theorem 
we first prove the following special case. Recall that
a ring is reduced if it has no non-zero nilpotent
elements.

\begin{thm} Let $R$ and $R^{'}$ be commutative rings and
assume that $R^{'}$ is reduced. If
$\phi:M_n(R)\rightarrow M_n(R^{'})$ is a ring homomorphism,
then $\phi(Z(M_n(R)))\subset Z(M_n(R^{'}))$.
\end{thm}

\begin{lem} Let $\phi:M_n(R)\rightarrow M_{n^{'}}(R^{'})$ be a
ring homomorphism. Then $n\leq n^{'}$.
\end{lem}

\noindent\textbf{Proof} i) Let $\p$ be a prime ideal in $R'$ and
denote by $\pi:M_{n'}(R')\rightarrow M_{n'}(R'/\p)$ the canonical
map sending each matrix $(a_{ij})$ to $(\overline{a_{ij}})$. 
Replacing $\f$ by $\pi\circ\f$ we may assume that $R'$
is a domain.\\ 
\noindent ii) We may also assume that $\phi$ is a
injective: namely let $J:=ker\phi$. Then $J$ is a 
two-sided ideal in the Azumaya algebra $A:=M_n(R)$, whence
by 2.3(2), $J=IA$ for some ideal $I$ in $R$.
So we get an induced injective ring homomorphism 
$M_n(R)/IM_n(R)\rightarrow M_{n^{'}}(R^{'})$. Since
 $M_n(R)/IM_n(R)\simeq M_n(R/I)$, we get an injective
 ring homomorphism from $M_n(R/I)\rightarrow M_{n^{'}}(R^{'})$.\\
\noindent ii) So let $\phi$ be injective and let $a\in M_n(R)$
be the standard Jordan cell of maximal rank i.e. the first
column of $a$ is zero and for each $i\geq 2$ the $i$-th column
of $a$  is equal to the $i-1$-th standard basis vector $e_{i-1}$.
So $a^n=0$ but $a^{n-1}\neq 0$. Consequently 
$\phi(a)^n=0$ and $\phi(a)^{n-1}\neq 0$ since $\phi$ is injective.
If now $n'<n$ it follows from the fact that $R'$ is a domain
 that $\phi(a)^{n-1}=0$, contradiction.

\medskip

\begin{rem}
\noindent The main interest of the above proof of 3.3 is that it is 
elementary. A shorter but less elementary proof consists in applying the
 converse of the Artin-Procesi theorem (2.6), using 2.2(1).
\end{rem}

\noindent\textbf{Proof of theorem 3.1}\\
\noindent i) First assume that $R^{'}$ is a domain. Then
$R^{'}\subset Q(R^{'})\subset k$, where
$k$ is an algebraic closure of the quotient field
$Q(R^{'})$ of $R^{'}$. Observe that $M_n(R^{'})\subset M_n(k)$
and that $M_n(R^{'})\cap Z(M_n(k))\subset Z(M_n(R^{'}))$. So we may assume
that $R^{'}=k$. Now let $c\in Z(M_n(R))$ and let
$\lambda$ be an eigenvalue of $\phi(c)$. Put

\medskip

$$V_{\lambda}:=\{v\in k^n| \phi(c)(v)=\lambda v\}.$$ 

\medskip

\noindent Since $c\in Z(M_n(R))$ we get that
 $\phi(c)\phi(x)=\phi(x)\phi(c)$ for all 
$x\in M_n(R)$, which implies that $\phi(x)V_{\lambda}
\subset V_{\lambda}$ for all $x\in M_n(R)$. So we get a ring
homomorphism

$$M_n(R)\ni x\rightarrow \phi(x)_{|V_{\lambda}}\in 
End(V_{\lambda})\simeq M_{n^{'}}(k)$$

\noindent where $n^{'}:=$dim$V_{\lambda}\leq n$. Since by
(3.2) $n\leq n^{'}$ we get $n=n^{'}$. So dim$V_{\lambda}=$
dim $k^{n}$ i.e. $V_{\lambda}=k^{n}$, whence $\phi(c)=\lambda I_n
\in Z(M_n(k))$.\\
\noindent ii) Now let $R^{'}$ be a reduced ring, $c\in Z(M_n(R))$
and $a\in M_n(R^{'})$. We need to show that $[\phi(c),a]=0$.
Therefore it suffices to show that for each $\p\in Spec(R^{'})$
all entries of the matrix $[\phi(c),a]$ belong to $\p$ (for
then these entries belong to $\cap\p=r(0)=(0)$ since
$R^{'}$ is reduced, so $[\phi(c),a]=0$). So let $\p\in Spec(R^{'})$
and denote by $\pi: R^{'}\rightarrow R^{'}/\p$ the canonical map.
Then extending $\pi$ to $M_n(R^{'})$ in the obvious way we get

$$\pi([\phi(c),a])=[(\pi\circ\phi)(c),\pi(a)]=0$$

\noindent by i), since $\pi\circ\phi:M_n(R)\rightarrow
M_n(R^{'}/\p)$ is a ring homomorphism, $c\in Z(M_n(R))$
and $R^{'}/\p$ is a domain. So all entries of $[\phi(c),a]$ belong to
$ker\pi=\p$, as desired.

\section{The main theorem}

The main result of this paper is the following theorem.

\begin{thm} Let $\f:A\rightarrow \apr$ be a ring homomorphism
between Azumaya algebras of constant ranks over
$R$ respectively $\rpr$.

\noindent (1) Then the rank of $A$ is lower or equal than the rank of $\apr$.

\noindent (2) If furthermore they are equal and  $\rpr$ is reduced, then 
$\f$ sends the center of $A$ to the center of $\apr$.
\end{thm}

\noindent\textbf{Proof} i) Let us first assume in addition that $R$ and $\rpr$ 
are 
finitely generated $\zx$-algebras, and let $\m\in Max(R^{'})$. According to 
2.3(2), there exists an ideal $I$ of $R$ such that ${\phi}^{-1}(\m 
A^{'})=IA$ . So we get an injective
map $\psi:A/IA\rightarrow A^{'}/\m A^{'}$ induced by $\f$ . Now observe 
that $A^{'}/\m A^{'}=A^{'}\tens_{R^{'}} R^{'}/\m$ is finite,
 since $R^{'}/\m$ is 
finite according to 2.10(2) and $A^{'}\tens_{R^{'}} R^{'}/\m R^{'}$ is 
finitely generated over $R^{'}/\m$ according to 2.3(1). Since $\psi$ is 
injective it follows that $A/IA$ is finite. By 2.3(4) both $A/IA$ and
$A^{'}/\m A^{'}$ are Azumaya algebras over $R/I$ respectively
$R^{'}/\m R^{'}$. So by 2.9 and the hypothesis on the rank of the 
Azumaya algebras $A$ and $A^{'}$ it follows that $A/IA\simeq M_n(R/I)$ and 
$A^{'}/\m A^{'}\simeq M_{n^{'}}(R^{'}/\m R^{'})$ for some integers $n\geq 1$ 
and $n^{'}\geq 1$. So statement (1) follows from 3.3.

\noindent ii) Furthermore by 3.1 we deduce that
$\psi(Z(A/IA))\subset Z(A^{'}/\m A^{'})$. Now let $c\in Z(A)$ and
$a^{'}\in A^{'}$. It follows that 
 $[\phi(c),a^{'}]\in\m A^{'}$ for all
$\m\in Max(R^{'})$. Hence by 2.3(3)

$$[\phi(c),a^{'}]\in\cap_{\m\in Max(R^{'})} \m A^{'}=
(\cap \m_{\m\in Max(R^{'})})A^{'}.$$

%\noindent By 2.2 (2) we have that $\m A^{'}\cap R^{'}=\m$, hence

%$$\cap_{\m\in Max(R^{'})} (\m A^{'}\cap R^{'})=\cap_{\m\in Max(R^{'})} \m$$

%So, again according to 2.2 (2), we get

%$$\cap_{\m\in Max(R^{'})}\m A^{'}=(\cap_{\m\in Max(R^{'})} \m A^{'} \cap 
%R^{'})A^{'}=
%(\cap_{\m\in Max(R^{'})} \m)A^{'}$$

\noindent Since $R^{'}$ is Jacobson by 2.10(2),
and reduced by hypothesis we have that

$$\cap_{\m\in Max(R^{'})} \m=\cap_{\p\in Spec(R^{'})}\p=(0)$$

\noindent It follows that $[\phi(c),a^{'}]=0$, as desired.\\
\noindent iii) In the general case of $R$ and $\rpr$ , let us consider $c\in 
Z(A)$ . According to Proposition 5.7, p.97 in [KO] there exist subrings 
$R_0$, $R_0^{'}$ of $R$ respectively $R^{'}$ which are finitely generated 
$\zx$-algebras and Azumaya algebras $A_0$ over $R_0$ and $A_0^{'}$ over 
$R_0^{'}$ such that $A\simeq A_0\otimes_{R_0} R$ and $A^{'}\simeq A_0^{'} 
\otimes_{R_0^{'}} R^{'}$.  So there is a subring $R_1$ of $R$ which 
contains $R_0$ and is  finitely generated over $\zx$ and there is a subring 
$R_1^{'}$ of $R^{'}$ which contains $R_0^{'}$ and is finitely generated 
over $\zx$ such that $c$ belongs to the image of the canonical embedding of 
$A_0\otimes_{R_0} R_1$ into $A$ and $\phi$ maps this image to the image of 
the canonical embedding of ${A_0}^{'}\otimes_{{R_0}^{'}} {R_1}^{'}$ into
 $A^{'}$. Since $A\simeq ({A_0}\otimes_{R_0} R_1)\otimes_{R_1} R$ and 
$A^{'}\simeq ({A_0}^{'}\otimes_{{R_0}^{'}} {R_1}^{'})\otimes _{{R_1}^{'}} R'$ 
statement (1) of 4.1 follows from i).\\
\noindent Finally it follows from ii) that $\phi(c)$ belongs to the canonical image in 
$A^{'}$ of the center of $A_0^{'}\otimes_{R_0^{'}} R_1^{'}$ . According to 
2.3(4), this means that $\phi(c)$ belongs to the center of $A^{'}$, as desired.

%\begin{rem}

\begin{rem} The main interest of the proof above is that it uses only tools of 
the classical theory of Azumaya algebras as exposed in [AG], [DI] and [KO].
The following proof is more direct and short, but less elementary. It does 
not need a reduction to finitely generated $\zx$-algebras and to finite 
rings. It consists in combining the converse of the Artin-Procesi theorem 
(2.6) with the eigenvalue argument in the proof of 3.1.
\end{rem}

\noindent\textbf{Sketch of an alternative proof of theorem 4.1}

\medskip

\noindent i) Statement (1) of the main theorem follows directly from the converse 
of the Artin-Procesi theorem since for any pair of positive integers $m$ and 
$n$ such that $m \leq n$, any element of $M_m(\zx)$ 
satisfies all $\zx$-identities of $M_n(\zx)$ (since $M_m(\zx)$ embeds
canonically into $M_n(\zx)$).\\
\noindent ii) Let $\f:A\rightarrow \apr $ be a ring morphism between Azumaya
 algebras of the same constant rank $n^2$ with center $R$ respectively $\rpr$
  such  that is $\rpr$ reduced.\\
\noindent  First assume that $R^{'}$ is a domain. Then
 $R^{'}\subset Q(R^{'})\subset k$, where $k$ is an algebraic closure of the
  quotient field $Q(R^{'})$ of 
$R^{'}$. Since the canonical map from $R'$ to $k$ is injective, it follows 
from the flatness of $A'$ over $R'$, (2.3(1)), 
that the canonical map $f$ from $A'=A^{'}\otimes_{R^{'}} R^{'}$ 
to $A^{'}\tens_{R^{'}} k$ is injective.
 On the other hand by 2.5(1) and 2.7(2), there exists an 
isomorphism $g$ of $k$-algebras from $A^{'}\tens_{R^{'}} k$ to 
$M_n(k)$. So if we put ${\f}' = g \circ f \circ \f$ we ghet a ring
homomorphism ${\f}':A\rightarrow M_n(k)$.
Furthermore, according to the injectivity of $g\circ f$ we get that
$A^{'}\cap Z(M_n(k))\subset Z(A')$. So we may assume that $A^{'}=M_n(k)$ and
$R^{'}=k$.
Now let $c\in R$ and  $\lambda$ be an eigenvalue of $\phi'(c)$. 
Following part i) of the proof of 3.1 (replacing everywhere $M_n(R)$
by $A$ and $\f$ by $\f^{'}$) we get a ring homomorphism
%Put
%$$V_{\lambda}:=\{v\in k^n| \phi'(c)(v)=\lambda v\}.$$
%By the definition of $R$, we get that $\phi'(c)\circ \phi'(x)=\phi'(x)\circ 
%\phi'(c)$ for any $x \in A$, which implies that $\phi'(x)V_{\lambda}
%\subset V_{\lambda}$ for any $x\in A$. So we get a ring
%homomorphism

$$A \ni x\rightarrow \phi'(x)_{|V_{\lambda}}\in
End(V_{\lambda})\simeq M_{n^{'}}(k)$$

\noindent where $n^{'}:=\dim V_{\lambda}\leq n$. According to  
statement (1) of 4.1 and 2.2(1), we have $n \leq n'$. So $n=n^{'}$
and as in the proof of 3.1 we deduce that $\phi'(c)=\lambda I_n
\in Z(M_n(k))$.\\
\noindent iii) Finally, let $R'$ be reduced, $c \in Z(A)$ and $a'\in A'$.
 We need to show that $[\phi(c),a']=0$.
So let $\p\in Spec(R^{'})$ and  $\pi:A^{'}\rightarrow A^{'}/\p A^{'}$ 
be the canonical map. Since by 2.3(4) $A^{'}/\p A^{'}$ is an Azumaya algebra
over $R^{'}/\p R^{'}$, which is a domain, it follows from ii) that
$\pi\circ\f:A\rightarrow A^{'}/\p A^{'}$ sends the center of $A$ to the center
of $A^{'}/\p A^{'}$. So $[\pi\circ\f (c),{\overline(a)}']=0$, for
all $\p\in Spec(R')$, whence by 2.3(3)

\medskip

$[\f (c),a^{'}]\in\cap (\p A^{'})=(\cap\p)A^{'}=(0)$

\medskip

\noindent since $R^{'}$ is reduced. So $[\phi(c),a']=0$ i.e. 
$\phi (c)\in Z(A^{'})$, as desired.

%\end{rem}

\begin{cor}
There is no algebra morphism from a Weyl algebra over a field $k$ to 
another Weyl algebra of strictly lower dimension over $k$.
\end{cor}

\noindent\textbf{Proof} i) If $k$ has positive characteristic
 it follows from 4.1 and 2.2(2).\\
\noindent ii) So assume that $k$ has characteristic zero and let  
 $\f:A_n(k)\rightarrow A_{n'}(k)$ be a morphism of $k$-algebras.
Then there exists a finitely generated $\zx$-algebra $R$ contained in $k$ 
such that $\f$ induces a morphism of $R$-algebra ${\f}_R:A_n(R)\rightarrow 
A_{n'}(R)$. Since $k$ has characteristic zero  
there exists a prime number $p$ whose canonical image in $R$ 
 belongs to a maximal ideal $\m'$ of $R$ (see for instance
 [ESS]. Prop.4.1.5). Let $k'$ be the residue field of $\m'$ and 
${\f}_{k'}:A_n(k')\rightarrow A_{n'}(k')$ the morphism of $k'$-algebras induced 
by ${\f}_R$. Since the characteristic of $k'$ is $p$, it follows from i) 
that $n \leq n'$ as desired.

\section{Isomorphisms of Azumaya algebras}

In this section we consider the following question:
when is a ring homomorphism $\phi$ between
Azumaya algebras an isomorphism ?\\
\noindent Apart from a technical assumption,
 the answer given below states that this
is the case if and only if $\phi$ induces an isomorphism
between the centers of the Azumaya algebras.
The proof is based on  the following easy generalization
 of a result given in [AG], whose proof we follow closely.

\begin{prop} Let $A$ and $A^{'}$ be Azumaya algebras over $R$
respectively $R^{'}$ and $\f:A\rightarrow A^{'}$ a ring homomorphism
such that the restriction of $\f$ to $R$ induces an 
isomorphism between $R$ and $R^{'}$. If 

\medskip

\noindent (*) rk$_{R_{\m}}A_{\m}$=rk$_{\rpr_{\f(\m)}}\apr_{\f(\m)} \mbox{ for all }  
\m\in Max(R)$

\medskip

\noindent then $\f$ is an isomorphism.
\end{prop} 

\noindent\textbf{Proof}. i) Observe that ker$\f$ is a two-sided ideal
in $A$. So by 2.3(2) it is of the form $IA$ for some ideal $I$ in $R$.
Let $i\in I$. Then $\f (i.1)=0$, so $i.1=0$ since $\f:R\rightarrow\rpr$
is injective. So $IA=0$ i.e. $\f$ is injective.\\
\noindent Consequently $\f:A\rightarrow A_2:=\f (A)$ is an isomorphism
and since $R'=\f (R)$, $A_2$ is an Azumaya algebra over $R'$ contained
in $A_1:=A'$ which is also an Azumaya algebra over $R'$. So by 2.2(8)
$A_1^{A_2}$ is an Azumaya algebra over $R'$ and
$\tau:A_2\otimes_{R'} A_1^{A_2}\rightarrow A_1$ is an isomorphism and by
2.2(1) $A_1^{A_2}=R'\oplus L$ for some $R'$-submodule of $L$. If we can show
that $L=0$ it follows that $A_1^{A_2}=R'$. Since $R'=\f (R)$ and $\tau$
is surjective this implies that $\f (A)=A'$.
%\noindent ii) It remains to see that $\f(A)=\apr$. Therefore consider

%$${\apr}^{\f (A)}:=\{b\in\apr |\, b\f (a)=\f (a)b \mbox{ for all } a \in A\}.$$

%\noindent It is shown in the proof of theorem 4.3 in [DI] that there
%exists an $\rpr$-submodule $L$ of $\apr$ such that $\apr=\f(A)\oplus L$.
%So it suffices to prove that $L_{\mpr}=0$ for all $\mpr\in Max(\rpr)$ (for
%then $L=0$ and hence $\apr=\f(A))$.
\noindent To see that $L=0$ let $\mpr\in Max(\rpr)$. Observe
 that by 2.3(1) $\apr$, $\f(A)$ and $A_1^{A_2}$ are finitely generated 
 projective $\rpr$-modules, and hence so is $L$. Consequently
the localizations of these modules with respect to
$\mpr$ are all free $\rpr_{\mpr}$-modules of finite rank, since $\rpr_{\mpr}$
is local. Now let $\m={\f}^{-1}(\mpr)$, so $\m\in Max(R)$. Since 
$\f:R\rightarrow\rpr$ is an isomorphism and $A\simeq\f(A)$,
 because $\f$ is injective,
we get that rk$_{R_{\m}}A_{\m}$=rk$_{\Rpr}\f(A)_{\mpr}$. Since 
by hypothesis rk$_{R_{\m}}A_{\m}$=rk$_{\Rpr}\apr_{\mpr}$ we get that
rk$_{\Rpr}\f(A)_{\mpr}$=rk$_{\Rpr}\apr_{\mpr}$. Then localizing the
isomorphism $\tau$ with respect to $\mpr$ and comparing the ranks of the
$\Rpr$-modules, we get that rk$_{{R'}_{\m '} }(A_1^{A_2})_{\m '}=1$.
Localizing the equation $A_1^{A_2}=R'\oplus L$ with respect to $\m '$
it follows that $L_{\m '}=0$ for all $\m '\in Max(R')$.
So $L=0$ as desired.

\medskip

\noindent Now we are able to give the main results of this section

\begin{cor} Any endomorphism of Azumaya algebras is an automorphism. 
\end{cor}

\begin{thm} A ring morphism $\f$ of an Azumaya algebra $A$ over $R$ to another Azumaya algebra
$A'$ over $R'$ is an isomorphism if and only if it induces an isomorphism between
$R$ and $R'$ and satisfies the condition (*) above.
\end{thm}

\noindent\textbf{Proof}. If $\f:A\rightarrow\apr$ is a ring isomorphism
then $\f$ is surjective. Hence it sends the center of $A$
to the center of $A^{'}$. So $\f$ induces a morphism from $R$ to
$\rpr$. Similarly ${\f}^{-1}$ induces an morphism from $\rpr$ to $R$.
 Consequently $\f$ induces an isomorphism from $R$ to $\rpr$ 
 and satisfies (*). Conversely, if $\f$ induces an isomorphism from $R$
  to $\rpr$ and satisfies (*) then it follows from 5.1
   that $\f:A\rightarrow\apr$ is an isomorphism.

\begin{cor} Let $\f:A\rightarrow\apr$ be a ring morphism between Azumaya
 algebras over $R$
 respectively $\rpr$ of the same constant rank. Then $\f$ is an isomorphism
  if and only if $\f$ induces an isomorphism between $R$ and $\rpr$.
\end{cor}

\noindent\textbf{Proof}. One easily verifies that our hypothesis on the
 ranks implies condition (*). So we can apply 5.3.

\section*{References}

\noindent [AG] M. Auslander and O. Goldman, The Brauer group of a commutative
ring, Trans. Amer. Math. Soc. 97 (1960), 367-409.

\medskip

\noindent [AR] M. Artin, On Azumaya Algebras and  Finite Dimensional Representations of 
Rings, Journal of Algebra 11 (1969), p. 532-563. 

\medskip

\noindent [AT] M. F. Attiya and I. G. Macdonald, Introduction to Commutative Algebra,
Addison-Wesley Publishing Compagny, 1969.

\medskip

\noindent [BO1] Bourbaki, Alg\`ebre Commutative, Masson, 1985.

\medskip

\noindent [BO2] Bourbaki, Alg\`ebre, Chap. VIII, Hermann, 1973.

\medskip

\noindent [DI] F. DeMeyer and E. Ingraham, Separable algebras over commutative rings,
Springer Lecture Notes 181, 1971.

\medskip

\noindent [DIE] J. Dieudonn\'e, De l'Analyse Fonctionnelle aux Fondements de la 
G\'eom\'etrie Alg\'ebrique, in the Grothendieck Festschrift, Vol. 1, Birkhauser, 1990, 
p. 1-14.

\medskip

\noindent [DIX] J. Dixmier, Sur les alg\` ebres de Weyl, Bull. Soc. Math. France 
96(1968), p. 209-242.

\medskip

\noindent [DR] V. Drensky and E. Formanek, Polynomial Identity Rings, 
Birkh\" auser-Verlag, 2004.

\medskip

\noindent [ESS] A. van den Essen, Polynomial Automorphisms and the
Jacobian Conjecture, Vol. 190 in Progress in Math., Birh\"auser, Basel, 2000.

\medskip

\noindent [KO] M-A. Knus and M. Ojanguren, Th\'eorie de la Descente et
Alg\`ebres d'Azumaya, Springer Lecture Notes 389, 1974.

\medskip

\noindent [M] J.C. MacConnell, J.C. Robson, Noncommutative 
Noetherian Rings, GMS 30, AMS, Providence, 2001.

\medskip

\noindent [O] M. Orzech and C. Small, The Brauer group of commutative rings, 
Lecture Notes in Pure and Applied Mathematics, Vol. 11, Marcel Dekker, 1975.

\medskip

\noindent [R] P. Revoy, Alg\`ebres de Weyl en caract\'eristique $p$, C.R. Acad.
Sci. Paris. S\'er. A-B 276 (1973), A 225-228.

\medskip

\noindent [RO] L. H. Rowen, Polynomial Identities in Ring Theory, Academic Press, 
New-York, 1980.

\medskip 

\noindent [P] C. Procesi, On a theorem of M. Artin, Journal of Algebra 22 (1972), 
p. 309-315.

\medskip

\noindent [T] Y. Tsuchimoto, Preliminaries on Dixmier Conjecture, Mem. Fac.
Sci. Kochi. Univ. Ser.A Math. 24 (2003), 43-59.

\medskip

\noindent [V] A. Verschoren, A check list on Brauer groups, in 
Brauer groups in ring theory and algebraic geometry (Wilrijk, 1981), p. 
279-300, Lecture Notes in Math. 917, Springer, 1982.

\medskip

\noindent Authors addresses

\medskip

\noindent Kossivi Adjamagbo\\
\noindent Universit\'e Paris 6\\
\noindent UFR 929, 4 Place Jussieu\\
\noindent 75252 Paris\\
\noindent France\\
\noindent Email: adja@math.jussieu.fr

\medskip

\noindent Jean-Yves Charbonnel\\
\noindent Universit\'e Paris 7-CNRS\\
\noindent Institut de Math\'ematiques de Jussieu\\
\noindent Th\' eorie des groupes\\
\noindent Case 7012, 2 Place Jussieu\\
\noindent 75251 Paris Cedex 05\\
\noindent France\\
\noindent Email: jyc@math.jussieu.fr

\medskip

\noindent Arno van den Essen\\
\noindent Department of Mathematics\\
\noindent Radboud University Nijmegen\\
\noindent Toernooiveld\\
\noindent 6525 ED Nijmegen\\
\noindent The Netherlands\\
\noindent Email: essen@math.ru.nl

\end{document}